\documentclass[journal]{IEEEtran}
\usepackage{bm}
\usepackage{amsmath,amsfonts,amsthm,amssymb}
\usepackage{booktabs}
\usepackage{footnote}
\usepackage[normalem]{ulem}
\makesavenoteenv{table}
\usepackage{url}
\usepackage{graphicx,float,wrapfig,subfigure}
\usepackage{xtab}
\usepackage{cite}
\usepackage{multirow}
\usepackage{multicol}
\usepackage{tikz}
\usepackage{algorithm, algorithmicx}
\usepackage[noend]{algpseudocode}
\usepackage{comment}
\makeatletter
\newcommand*{\addFilecorrelation}[1]{% argument=file name and extension
	\typeout{(#1)}% latexmk will find this if $recorder=0 (however, in that case, it will ignore #1 if it is a .aux or .pdf file etc and it exists! if it doesn't exist, it will appear in the list of dependents regardless)
	\@addtofilelist{#1}% if you want it to appear in \listfiles, not really necessary and latexmk doesn't use this
	\IfFileExists{#1}{}{\typeout{No file #1.}}% latexmk will find this message if #1 doesn't exist (yet)
}
\makeatother

\theoremstyle{plain}

\theoremstyle{remark}

\theoremstyle{remark}
\newtheorem{example}{Example}

\graphicspath{{figures/}}

\includecomment{extended}
\excludecomment{reduced}

\begin{document}

\title{
    Optimal Control of DERs in ADN under Spatial and Temporal Correlated Uncertainties
	}

\author{
Xiaoshuang~Chen,~\IEEEmembership{Student Member,~IEEE,}
Jin~Lin,~\IEEEmembership{Member,~IEEE,}
Feng~Liu,~\IEEEmembership{Member,~IEEE,}
and~Yonghua~Song,~\IEEEmembership{Fellow,~IEEE}

\thanks{
%	This work was supported by National High-Technology Research and Development Program (“863” Program) of China (2014AA051901), National Natural Science Foundation of China (51577096,51761135015), International S\&T Cooperation (2016YFE0102600).
	
    X. Chen, J. Lin, F. Liu and Y. Song are with the State Key Laboratory of Control and Simulation of Power Systems and Generation Equipment, Department of Electrical Engineering, Tsinghua University, Beijing 100084, China. (email: linjin@tsinghua.edu.cn).

%    Y. Song is with the Department of Electrical and Computer Engineering, University of Macau, Macau, China, and the Department of Electrical Engineering, Tsinghua University, Beijing 100084, China.
	}
}

\maketitle
\begin{abstract} 
    The control schemes of distributed energy resources (DERs) in active distribution networks (ADNs) are largely influenced by uncertainties. The uncertainties of DERs are complicated, containing spatial and temporal correlation, which makes it challenging to design proper control schemes, especially when there exist temporal-correlated units such as energy units (EUs). This paper provides an It\^{o} process model to describe the characteristics of stochastic resources and EUs in a unified way, which makes it easy to evaluate the impacts of stochastic resources on temporal-correlated units. Based the moment form of the It\^{o} process model, a moment optimization (MO) approach is provided to transform the stochastic control (SC) problem into an optimization problem with respect to the first-order and second-order moments of the system variables. The scale of MO is comparable to that of the corresponding deterministic control problem, which means that the computational efficiency of MO is much smaller than that of traditional approaches.  Case studies also show that the proposed approach outperforms existing approaches in both the performance and computational efficiency, which means that the proposed approach has attractive potential for use in large-scale applications.
\end{abstract}
%\vspace{-3mm}

\vspace{-2mm}
\section{Introduction}
\IEEEPARstart{D}{istributed} energy resources (DERs), including renewable generations and energy units (EUs), have grown rapidly in recent years \cite{maycock2014future}. The integration of DERs brings significant challenges to active distribution networks (ADNs), such as overvoltage and overloading problems \cite{hidalgo2010review}. To this end, it is necessary to control DERs in order to mitigate their negative impacts, and there are many studies in this area \cite{olivier2016active,li2016convex,sugihara2013economic}. However, renewable power generations are usually highly stochastic, and the uncertainties caused by these renewables will undoubtedly deteriorate the control performance.

It is challenging to consider DER uncertainties in the optimal control problem, mainly because the modeling of these uncertainties is complicated. The uncertainty of renewable generations is usually non-Gaussian \cite{bludszuweit2008statistical, li2015data,teng2013optimal,atwa2010optimal} and contains spatial and temporal correlations \cite{tabone2015modeling,lorenz2009irradiance}.
Moreover, EU characteristics are also temporal-correlated, and their temporal correlations may be influenced by the spatial and temporal correlation of the uncertainties \cite{chen2018unified}. However, based on existing uncertainty models, such as the probability distribution model \cite{bludszuweit2008statistical} and the Markov model \cite{tabone2015modeling}, it is challenging to analyze the impacts of uncertainties on the DERs and distribution networks, and such impacts can only be evaluated by Monte Carlo simulations \cite{anese2015optimal}, which are time-consuming.

Techniques to deal with uncertainties in optimal control problems include robust control methods\cite{li2018coordinated}, model predictive control (MPC) methods \cite{valverde2013model,xing2017model,meng2015cooperation}, and stochastic-programming-based-control (SPBC) methods \cite{anese2015optimal,atwa2010optimal,agalgaonkar2015stochastic}. Robust controllers use uncertainty sets to model uncertainties and find the control schemes that perform well in the worst case; hence, these schemes usually lead to conservative results. MPC solves an open-loop optimal control problem in which the uncertainties are not considered and adjusts the control outputs in a receding-horizon manner. Although it provides some robustness by the receding-horizon implementation, MPC does not consider the uncertainties explicitly, which may have negative performance impacts \cite{mesbah2016stochastic}. Moreover, receding-horizon implementation is also time-consuming.

SPBC has been widely used in the recent studies \cite{anese2015optimal,atwa2010optimal,agalgaonkar2015stochastic}; SPBC handles the stochastic control (SC) problem by stochastic programming, which can be solved by scenario-based approaches. Specifically, SPBC generates a certain number of scenarios under the probability distribution and correlation of the uncertainties; furthermore, it transforms the SC into deterministic optimization problems. SPBC is widely used in power system operations considering uncertainties; however, in order to achieve good accuracy, a large number of scenarios are needed, which may lead to an unacceptable computational burden. Although some studies exist for methods to accelerate the computation of SPBC \cite{fu2016multiobjective,zhu2014decomposition}, the computational burden of SPBC is still too large compared to deterministic control.

In summary, it is challenging to efficiently solve the SC problem under complicated uncertainties with spatial and temporal correlation. Therefore, this paper provides a novel moment optimization (MO) approach for the SC of DERs in distribution networks. We use an It\^{o} process model to describe the probability distribution, the spatial and temporal correlation of uncertainties, and transform the SC into a deterministic optimization with respect to the first-order and second-order moments of system variables. The proposed MO approach solves the SC with a comparable computational burden to deterministic control problems and hence is attractive for online applications. Case studies also show that the MO approach outperforms the existing approaches. 

The contributions of this paper are twofold:
\begin{enumerate}
	\item An It\^{o} process model is provided to describe the stochastic resources. On the one hand, the proposed model can be used to model the spatial and temporal correlation of renewable generations; on the other hand, the temporal correlation of EUs can be easily embedded into the It\^{o} process model; therefore, it is possible to consider the impact of renewable generations and control policies on the states of EUs in a unified framework. Moreover, the statistics of It\^{o} processes can be calculated analytically, i.e., without time-consuming simulations.
	\item An SC model of DERs is provided and then solved by the MO approach. The MO approach transforms the SC into a deterministic optimization problem with respect to the first-order and second-order moments of the system variables. The MO model accurately describes the characteristics of the system, and its scale is comparable to that of the corresponding deterministic control problems. Therefore, MO achieves a good tradeoff between performance and computational efficiency.
\end{enumerate}

Following this introduction, Section \ref{section:model} provides the model of the stochastic resources and the SC problem. The MO approach is discussed in detail in Section \ref{section:mo}. Section \ref{section:case} provides numerical results, and Section \ref{section:conclusion} concludes the paper.
\begin{reduced}
	Due to space constraints, some proofs are put into an extended version which is accessible online \cite{Chen2019Optimal}.
\end{reduced}

\vspace{-2mm}
\section{Modeling} \label{section:model}
This section describes the SC model of distribution networks. Typical units as well as the uncertainties in distribution networks are considered. After discussing the structure of ADN, we provide the It\^{o} process model of stochastic resources and then establish the SC model. Although we use a continuous-time formulation for convenience, all models can be easily transformed into a discrete-time formulation.

\vspace{-2mm}
\subsection{Brief Structure of ADN}
\begin{figure}
	\centering
	\includegraphics[width=\columnwidth]{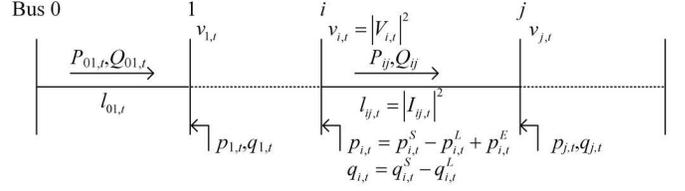}
	\caption{Brief structure of a radial ADN.}
	\label{fig:adn-structure}
	\vspace{-4mm}
\end{figure}

Fig. \ref{fig:adn-structure} shows the brief structure of a radial ADN. The tree topology of ADN is described by a set of buses, denoted by $\mathcal{V} = \left\{0,1,\cdots,N\right\}$, and a set of branches, denoted by $\mathcal{E}=\left\{(i,j)\right\}$. Moreover, a set of time is denoted by $\mathcal{T}$. In general, we regard $i,j,k$ as the bus indices and $t$ as the time index. Moreover, we set Bus 0 as the root bus connected to the external grid.

For Bus $i$, denoted by $v_{i,t}$ the square voltage amplitude, and $p_{i,t}$ and $q_{i,t}$ are the active and reactive power injection, respectively. For Branch $(i,j)$, $P_{ij,t}$ and $Q_{ij,t}$ denote the active and reactive power flow, respectively, from Bus $i$ to Bus $j$, and $l_{ij,t}$ denotes the square current amplitude.

There are three types of power injections:
\begin{enumerate}
	\item Fixed load, of which the active and reactive power are denoted by $p_{i,t}^L$ and $q_{i,t}^L$, respectively. $p_{i,t}^L$ and $q_{i,t}^L$ are considered uncontrollable and are assumed not to have uncertainties. The traditional load at each bus is regarded as this type.
	\item Stochastic resources such as renewable generations. The active power of the stochastic resources, denoted by $p_{i,t}^S$, contains uncertainty. However, the reactive power, denoted by $q_{i,t}^S$, is considered controllable in ADN due to the fact that the grid-connected converter is able to adjust the reactive power output. 
	\item EUs including battery storage and thermostatically controlled loads. The active power of the $i$-th EU is denoted by $p_{i,t}^E$. An important constraint of EUs is the state-of-charge (SOC), denoted by $SOC_{i,t}$.
\end{enumerate}

The following subsections will provide the model of the ADN and the control problem in detail.

\vspace{-2mm}
\subsection{It\^{o} Process Model of Stochastic Resources} \label{section:ito}
Here we consider the model of stochastic resource $p_{i,t}^S$. The active power of stochastic resources can be separated into the prediction and deviation parts, denoted by $p^{pred}_{i,t}$ and $p^{dev}_{i,t}$, respectively. We regard the prediction part as fixed values in the SC problem. For convenience, let $\bm{\xi}_t = \left(p_{i,t}^{dev}\right)_{i\in\mathcal{V}}$.

Different from existing studies \cite{bludszuweit2008statistical,tabone2015modeling}, an It\^{o} process model is used in this work to describe the characteristics of $\bm{\xi}_t$. On the one hand, the model can describe the spatial and temporal correlation, as well as the probability distribution of renewables; on the other hand, the It\^{o} process, described by an SDE, is compatible with the description of EUs (see \eqref{eq:eu-dynamic}), which makes it possible to embed the stochastic characteristics into the optimization problem without Monte Carlo simulations. 
% Quality control editor: Abbreviations and acronyms (e.g., SDE) are typically defined the first time the term is used within the main text and then used throughout the remainder of the manuscript. Please consider adhering to this convention. The target journal may have a list of abbreviations that are considered common enough that they do not need to be defined.
In contrast, existing models need to be broken into a number of scenarios to be used in SC problems, which leads to unbearable computational burden.

The It\^{o} process model is defined as
\vspace{-0.5mm}
\begin{align}
%\bm{z}_t &= \bm{C}\bm{\xi}_t \label{eq:z-xi-vector} \\
d\bm{\xi}_t &= \bm{\mu}(\bm{\xi}_t)dt + \bm{\sigma}(\bm{\xi}_t)d\bm{W}_t \label{eq:xi-ito-vector}
\end{align}
where $\bm{\mu}\left(\cdot\right)$ and $\bm{\sigma}\left(\cdot\right)$ are the drift function and the diffusion function, respectively. Note that \eqref{eq:xi-ito-vector} is a stochastic differential equation (SDE); hence, we also need the initial condition describing the distribution of $\bm{\xi}_0$ in order to fully describe the It\^{o} process. Here, we omit it for convenience.

By properly setting $\bm{\mu}(\cdot)$ and $\bm{\sigma}(\cdot)$, the It\^{o} process model describes a large class of stochastic processes \cite{lamberton2011introduction} with different probability distribution and spatial/temporal correlation, which is also supported by the authors' previous work \cite{chen2018stochastic}. Here we provide some examples. 

\begin{example} \label{eg:ito-1}
	Consider the following It\^{o} process:
	\vspace{-0.5mm}
	\begin{align} \label{eq:ito-example-1}
	d\bm{\xi}_t &= -\frac{1}{\tau}\bm{\xi}_tdt + \frac{1}{\sqrt{\tau}}\bm{\sigma}d\bm{W}_t
	\end{align}
	i.e., $\bm{\mu}\left(\bm{\xi}_t\right)=-\bm{\xi}_t/\tau,\bm{\sigma}\left(\bm{\xi}_t\right)=\bm{\sigma}/\sqrt{\tau}$. Assume $\tau>0,\bm{\xi}_0=\bm{0}$.
	
	For this linear SDE, its second-order moment matrix is:
	\begin{equation}
	\mathbb{E}\bm{\xi}_t\bm{\xi}_s^\top = \frac{\exp\left(-\frac{t-s}{\tau}\right)-\exp\left(-\frac{t+s}{\tau}\right)}{2}\bm{\sigma\sigma}^\top,\forall s\leq t \label{eq:xi-ito-example}
	\end{equation}
	It is easy to conclude from \eqref{eq:xi-ito-example} the spatial and temporal correlations of $\bm{\xi}_t$:
	\begin{itemize}
		\item[-] Spatial correlation: In practice, we usually use $\mathbb{E}\bm{\xi}_t\bm{\xi}_t^\top$ to describe the spatial correlation \cite{tabone2015modeling}, which is $\left[1-\exp\left(-2t/\tau\right)\right]\bm{\sigma\sigma}^\top/2$. When $t$ is large, the covariance matrix is approximately $\bm{\sigma\sigma}^\top/2$.
		\item[-] Temporal correlation: the temporal correlation is proportional to the term $\exp\left(-(t-s)/\tau\right)-\exp\left(-(t+s)/\tau\right)$. Specifically, when $t,s$ are sufficiently large, $\exp\left(-(t+s)/\tau\right)$ can be ignored, and the covariance is only determined by $t-s$, which means that $\bm{\xi}_t$ is a stationary process \cite{chen2018unified}. It is clear that the larger $\tau$ is, the larger the temporal correlation is.
	\end{itemize}
\end{example}

In summary, in this simple example, $\tau$ describes the temporal correlation, and $\bm{\sigma}$ describes the spatial correlation.

\begin{example}
	Consider the following It\^{o} process:
	\begin{align} \label{eq:ito-example-2}
	d\bm{\xi}_t &= -\frac{1}{\tau}\left(\bm{\xi}_t-0.5\right)dt + \frac{1}{\sqrt{\tau}}\sqrt{\bm{\xi}_t\left(1-\bm{\xi}_t\right)^\top}\bm{\sigma}d\bm{W}_t
	\end{align}
	This nonlinear SDE is similar to \eqref{eq:ito-example-1}, but it can be proved that each element of $\bm{\xi}_t$ follows a Beta distribution rather than a Gaussian distribution \cite{chen2018stochastic}. Actually, proper formulations of $\bm{\mu}(\cdot)$ and $\bm{\sigma}\left(\cdot\right)$ can describe a large class of probability distributions, such as Gaussian distribution, Beta distribution, Laplace distribution, and Weibull distribution \cite{chen2018stochastic}.
\end{example}

\begin{example}
	Here we consider a general form
	\begin{align}
	d\bm{\xi}_t &= \bm{\mu}(\bm{\xi}_t;\bm{\Theta})dt + \bm{\sigma}(\bm{\xi}_t;\bm{\Theta})d\bm{W}_t \label{eq:ito-example-3}
	\end{align}
	where $\bm{\mu}(\bm{\xi}_t;\bm{\Theta})$ and $\bm{\sigma}(\bm{\xi}_t;\bm{\Theta})$ are the parametrized drift function and diffusion function, respectively, and $\bm{\Theta}$ is the parameter to be determined. $\bm{\mu}$ and $\bm{\sigma}$ can, of course, be nonlinear, hence \eqref{eq:ito-example-3} can be used to describe general stochastic resources. Moreover, there are a number of studies on statistically determining $\bm{\Theta}$ based on historical data via maximum likelihood estimation \cite{du2017parameter,ait2002maximum,ait2008closed,yoshida1992estimation}. The computation of the spatial and temporal correlations of \eqref{eq:ito-example-3} is provided in \cite{chen2018stochastic}, which is further discussed in Section \ref{section:mo-stochastic}.
\end{example}

  In the remainder of this paper, we regard $\bm{\mu}(\cdot)$ and $\bm{\sigma}\left(\cdot\right)$  as a priori knowledge obtained by historical data.

\subsection{Stochastic Control Problem}
The SC problem aims at decreasing both the cost of electricity and the voltage profile under the constraints of state variables and controllable units. Specifically, the SC problem can be formulated as \eqref{eq:scp-obj}$\sim$\eqref{eq:sc-inequality}:
\begin{align} \label{eq:scp-obj}
\begin{aligned}
\min_{\bm{u}_t:t\in\mathcal{T}} J =& \mathbb{E}_{\bm{\xi}_t}\left\{\int_{t\in\mathcal{T}} \lambda_tP_{01,t}dt +R^V\sum_i\int_{t\in\mathcal{T}}(v_{i,t}-1)^2dt\right\} \\
&+ \mathbb{E}_{\bm{\xi}_t}\left\{\int_{t\in\mathcal{T}}\bm{u}_t^\top \bm{R}^U \bm{u}_tdt+R^E\sum_iSOC_{i,T}^2\right\}
\end{aligned}
\end{align}

Stochastic Resources and Control Variables:
\begin{subequations}
\begin{flalign}
& p^S_{i,t} = p^{pred}_{i,t} + p^{dev}_{i,t}, \forall i\in \mathcal{V}, t\in\mathcal{T}&& \label{eq:scp-p-s}\\
&\bm{\xi}_t = \left(p_{i,t}^{dev}\right)_{i\in\mathcal{V}},\forall t\in\mathcal{T}&& \label{eq:scp-vector-z-t}\\
&\bm{u}_t = \left(p_{i,t}^E,q_{i,t}^S\right)_{i\in\mathcal{V}},\forall t\in\mathcal{T}&&\label{eq:scp-vector-u-t}
\end{flalign}
\end{subequations}

EUs:
%\begin{subequations}
\begin{flalign}
&\frac{d}{dt}{SOC}_{i,t} = -\alpha_i^ESOC_{i,t}^E + \beta_i^Ep_{i,t}^E,\forall i\in\mathcal{V},t\in\mathcal{T}&&\label{eq:eu-dynamic}
\end{flalign}

Power and Network Constraints:
%\vspace{-1mm}
\begin{subequations} \label{eq:sc-network}
\begin{flalign}
&p_{i,t} = p_{i,t}^S-p_{i,t}^L + p_{i,t}^E, \forall i\in \mathcal{V},t\in\mathcal{T}&& \label{eq:scp-p-i-t}\\
&q_{i,t} = q_{i,t}^S - q_{i,t}^L, \forall i\in \mathcal{V}, t\in\mathcal{T}&& \label{eq:scp-q-i-t} \\
&p_{i,t}=\sum_{j:i\to j}P_{ij,t} -\sum_{k:k\to i} \left(P_{ki,t}-r_{ki}l_{ki,t}\right) + g_iv_{i,t},\forall i \in \mathcal{V}&& \label{eq:network-p} \\
&q_{i,t}=\sum_{j:i\to j}Q_{ij,t} -\sum_{k:k\to i} \left(Q_{ki,t}-x_{ki}l_{ki,t}\right) + b_iv_{i,t},\forall i \in \mathcal{V}&& \label{eq:network-q} \\
&\begin{aligned}
v_{j,t}= v_{i,t} &- 2\left(r_{ij}P_{ij,t}+x_{ij}Q_{ij,t}\right) \\&+ \left(r_{ij}^2+x_{ij}^2\right)l_{ij,t},\forall (i,j)\in \mathcal{E}
\end{aligned}&& \label{eq:network-v}\\
&l_{ij,t}v_{i,t}=P_{ij,t}^2+Q_{ij,t}^2,\forall (i,j)\in \mathcal{E}&& \label{eq:network-lv-pq} 
\end{flalign}
\end{subequations}

%\end{subequations}
%\vspace{-1mm}
Inequality Constraints (with the confidence level $\gamma$):
\begin{subequations} \label{eq:sc-inequality}
\begin{flalign}
&\left(p_{i,t}^S\right)^2 + \left(q_{i,t}^S\right)^2 \leq \left(\bar{s}_{i,t}^S\right)^2,\forall i\in\mathcal{V},t\in\mathcal{T}&&\label{eq:p-s-limit-original}\\
&\underline{SOC}_{i}\leq SOC_{i,t} \leq \overline{SOC}_{i}\forall i\in\mathcal{V},t\in\mathcal{T}&& \label{eq:eu-limit-e}\\
&\underline{p}_{i}^E\leq p_{i,t}^E \leq \bar{p}_{i}^E\forall i\in\mathcal{V},t\in\mathcal{T}&& \label{eq:eu-limit-p}\\
&\underline{v}_i\leq v_{i,t} \leq \bar{v}_i,\forall i\in\mathcal{V},t\in\mathcal{T}&& \label{eq:voltage-limit}\\
&0 \leq l_{ij,t} \leq \bar{l}_{ij},\forall (i,j)\in\mathcal{E},t\in\mathcal{T}&& \label{eq:current-limit}
\end{flalign}
\end{subequations}

%\vspace{-1mm}
In this paper, $\mathcal{T}$ is regarded as a finite interval; hence, the abovementioned model is a finite-horizon stochastic control model. The decision variable is $\bm{u}_t$, including the reactive power of stochastic resources and the active power of EUs. The objective is in the sense of expectation, and the constraints are classified into the following two groups: the equality constraints, including \eqref{eq:eu-dynamic}\eqref{eq:sc-network}, are regarded as almost sure constraints; while the inequality constraints in \eqref{eq:p-s-limit-original}$\sim$\eqref{eq:current-limit} are regarded as chance constraints, i.e., each inequality is satisfied with the confidence level $\gamma$.

These equations are explained in detail below. 

\subsubsection{Objective}
The objective \eqref{eq:scp-obj} of the SC problem is to minimize the expected cost function, which contains 4 parts:
%\vspace{-1mm}
\begin{itemize}
	\item[-] The price of electricity bought from the market, i.e., $\int_{t\in\mathcal{T}} \lambda_tP_{01,t}dt$, where $\lambda_t$ is the electricity price.
	\item[-] The penalty of the voltage profile $\sum_i\int_{t\in\mathcal{T}}(v_{i,t}-1)^2dt$.
	\item[-] The cost of control $\int_{t\in\mathcal{T}}\bm{u}_t^\top \bm{R} \bm{u}_tdt$.
	\item[-] The penalty of SOC $R^E\sum_iSOC_{i,T}^2$. This term drives the EUs to the default state in order to maximize the ability of operation in the future.
\end{itemize}
%\vspace{-1mm}

\subsubsection{Energy Units}
EUs are units with energy constraints, including battery storage systems, demand response, etc. 
% Quality control editor: Abbreviations and acronyms typically need to be defined only once within the main text. Please consider adhering to this convention.
\eqref{eq:eu-dynamic} describes the dynamics of EUs, where $\alpha_i^E$ is the dissipation factor, and $\beta_i^E$ is the charging efficiency.

\subsubsection{Power and Network Constraints}
\eqref{eq:scp-p-i-t} and \eqref{eq:scp-q-i-t} are the power balance equations at each bus. Here we adopt the distFlow model \cite{farivar2013branch} of network constraints, as shown in \eqref{eq:network-p} $\sim$ \eqref{eq:network-lv-pq}, where $r_{ij}$ and $x_{ij}$ are the resistance and reactance of the branch from Bus $i$ to Bus $j$, respectively; $g_i$ and $b_i$ are the shunt conductance and susceptance at Bus $i$, respectively.

\subsubsection{Inequality Constraints}
\eqref{eq:p-s-limit-original} is the capacity limit of the renewable generations. According to \cite{xing2017model}, this convex constraint can be approximately described by polygons, i.e.,
\begin{equation} 
\bm{C}^S_{i}\left[p_{i,t}^S, q_{i,t}^S\right]^\top\leq \bm{D}_i^S, \forall i \in\mathcal{V},t\in\mathcal{T} \label{eq:p-s-limit}
\end{equation}

\eqref{eq:eu-limit-e} and \eqref{eq:eu-limit-p} are the energy and power constraints of the EUs, respectively; \eqref{eq:voltage-limit} and \eqref{eq:current-limit} are the voltage and current constraints, respectively.

\vspace{-2mm}
\subsection{Compact Reformulation}
For convenience, we use a group of vectors to represent the abovementioned variables and transform the equations into a compact form. 

\subsubsection{Groups of Variables}
We have defined $\bm{\xi}_t$ and $\bm{u}_t$ in \eqref{eq:scp-vector-z-t} and \eqref{eq:scp-vector-u-t}. Now we define the following vectors, all of which are formulated as column vectors:
\begin{equation}
\begin{split}
\bm{e}_t &= \left(SOC_{i,t}^E\right)_{i\in\mathcal{V}} \\
\bm{x}_t &= \left(\left.P_{ij,t}\right|_{(i,j)\in\mathcal{E}},\left.Q_{ij,t}\right|_{(i,j)\in\mathcal{E}},\left.v_{i,t}\right|_{i\in\mathcal{V}}\right) \\
\bm{y}_t &= \left(l_{ij,t}\right)_{(i,j)\in\mathcal{E}} \\
\bm{d}_t &= \left(p_{i,t}^{pred},p_{i,t}^L,q_{i,t}^C\right)_{i\in\mathcal{V}}
\end{split}
\end{equation}
where $\bm{e}_t$ is the vector of the energy of EUs; $\bm{x}_t$ is the vector of network states except the branch currents, and $\bm{y}_t$ is the vector of branch currents; and $\bm{d}_t$ is the vector of fixed values, including prediction of renewables, load profile, and the reactive power supply of each bus.

\subsubsection{Affine Feedback Control Policy}
There are two kinds of control policies: state-feedback control policies and disturbance feedback control policies. Under certain circumstances, they are equivalent \cite{skaf2010design}. Here, we adopt the affine disturbance feedback control policy, i.e., $\bm{u}_t = \bm{u}_t^0 + \bm{K}\bm{\xi}_t$, where $\bm{u}_t^0$ and $\bm{K}$ are the decision variables.

\subsubsection{Compact Form of SC}
By the abovementioned vectors, it is possible to reformulate the SC model as
\begin{subequations} \label{eq:sc-vector}
\begin{equation} \label{eq:scp-vector-obj}
\begin{split}
\min_{\bm{u}_t^0,\bm{K}} J =& \mathbb{E}_{\bm{\xi}_t}\left\{\int_{t\in\mathcal{T}}\left(\bm{H} \bm{x}_t + \bm{x}_t^\top\bm{R}^X\bm{x}_t\right)dt\right\} \\
&+\mathbb{E}_{\bm{\xi}_t}\left\{\int_{t\in\mathcal{T}} \bm{u}_t^\top \bm{R}^U \bm{u}_tdt+\bm{e}_T^\top\bm{R}^E\bm{e}_T\right\}
\end{split}
\end{equation}
\begin{align}
\text{s.t.~}&d\bm{\xi}_t = \bm{\mu}(\bm{\xi}_t)dt + \bm{\sigma}(\bm{\xi}_t)d\bm{W}_t \label{eq:sc-vector-ito}\\
&\bm{u}_t = \bm{u}_t^0 + \bm{K}\bm{\xi}_t\label{eq:sc-control-policy}\\
&\dot{\bm{e}}_t = -\bm{\alpha e}_t + \bm{\beta u}_t\label{eq:eu-vector} \\
&\bm{x}_t = \bm{A}_y\bm{y}_t +\bm{A}_\xi\bm{\xi}_t+\bm{A}_d\bm{d}_t + \bm{A}_e\bm{e}_t + 
\bm{A}_u\bm{u}_t\label{eq:network-vector} \\
&l_{ij,t}v_{i,t}=P_{ij,t}^2+Q_{ij,t}^2,\forall (i,j)\in \mathcal{E} \label{eq:network-lv-pq-2} \\
&\underline{\bm{x}} \leq \bm{x}_t \leq \bar{\bm{x}} \label{eq:sc-constraints-x} \\
&\underline{\bm{y}} \leq \bm{y}_t \leq \bar{\bm{y}} \label{eq:sc-constraints-y}\\
&\underline{\bm{e}} \leq \bm{e}_t \leq \bar{\bm{e}} \label{eq:sc-constraints-e}\\
&\bm{C}_u\bm{u}_t + \bm{C}_\xi\bm{\xi}_t \leq \bm{D} \label{eq:sc-constraints-p-s}
\end{align}
\end{subequations}
where \eqref{eq:scp-vector-obj} is corresponding to \eqref{eq:scp-obj}; \eqref{eq:sc-control-policy} is the disturbance feedback control policy; \eqref{eq:eu-vector} corresponds to \eqref{eq:eu-dynamic};  \eqref{eq:network-vector} corresponds to \eqref{eq:sc-network}; and \eqref{eq:sc-constraints-x}$\sim$\eqref{eq:sc-constraints-p-s} correspond to \eqref{eq:sc-inequality}.

Here we provide more explanations to \eqref{eq:network-vector} and \eqref{eq:network-lv-pq-2}. Note that \eqref{eq:network-lv-pq-2} is the same as \eqref{eq:network-lv-pq}, which is not in a vector form. It is because \eqref{eq:network-lv-pq} is the only nonlinear constraint, which we address separately in Section \ref{section:mo}. To obtain \eqref{eq:network-vector}, by substituting \eqref{eq:scp-p-s}\eqref{eq:scp-p-i-t}\eqref{eq:scp-q-i-t} into \eqref{eq:network-p}\eqref{eq:network-q}\eqref{eq:network-v}, using the vector notations, we have
\begin{equation} \label{eq:network-vector-origin}
\bm{A}_x^0\bm{x}_t +\bm{A}_y^0\bm{y}_t +\bm{A}_\xi^0\bm{\xi}_t+\bm{A}_d^0\bm{d}_t + \bm{A}_e^0\bm{e}_t + \bm{A}_u^0\bm{u}_t = 0
\end{equation}
where it is easy to find that $\bm{A}_x^0$ is invertible. Therefore, when multiplied by $\left(\bm{A}_x^0\right)^{-1}$ in \eqref{eq:network-vector-origin}, we obtain \eqref{eq:network-vector}, where $\bm{A}_y^0 = -\left(\bm{A}_x^0\right)^{-1}\bm{A}_y^0$, and $\bm{A}_{\xi},\bm{A}_d,\bm{A}_e,\bm{A}_u$ are similarly obtained.

\eqref{eq:sc-vector} is the control model used for analysis in this paper. Similar models are also studied in \cite{valverde2013model,xing2017model,meng2015cooperation}, except for the model of stochastic resources \eqref{eq:sc-vector-ito}. Existing approaches usually use SPBC to solve the optimal control problem with stochastic resources. However, when considering spatial and temporal correlations, it is necessary to use a large number of scenarios to guarantee the accuracy, which makes the optimization undoubtedly time-consuming. Based on this, Section \ref{section:mo} provides the MO approach to efficiently solve \eqref{eq:sc-vector}.

\vspace{-2mm}
\section{Solution Based on Moment Optimization} \label{section:mo}
This section provides the MO approach to solving the SC problem. The basic idea of MO is based on the fact that the objective in \eqref{eq:scp-vector-obj} only contains quadratic forms. Therefore, the objective can be equivalently transformed into the function of the first-order and second-order moments of these variables. Furthermore, we can regard these first-order and second-order moments as decision variables, and all we have to do is find the the constraints of these moments. In other words, higher-order moments will not influence the solution of SC. By transforming the SC problem \eqref{eq:sc-vector} into a deterministic optimization problem with respect to the first-order and second-order moments, the MO approach largely reduces the computational burden of SC compared to traditional scenario-based approaches, without sacrificing performance.

We define the first two notations. Assume that $a$ is a certain component of $\bm{x}_t,\bm{y}_t,\bm{z}_t,\bm{u}_t,\bm{e}_t$ and denote by $\tilde{a}$ the expectation of $a$ and $\hat{a}$ the standard deviation of $a$, i.e.,
\begin{equation}
\begin{split}
\tilde{a} &= \mathbb{E}\left\{a\right\} \\
\hat{a} &= \sqrt{\text{var}(a)} = \sqrt{\mathbb{E}\left\{a^2\right\}-\left(\mathbb{E}\left\{a\right\}\right)^2}
\end{split}
\end{equation}
With these notations, we can define the first-order moments $\tilde{\bm{x}}_t,\tilde{\bm{y}}_t,\tilde{\bm{z}}_t,\tilde{\bm{u}}_t,\tilde{\bm{e}}_t$, and the second-order central moments $\hat{\bm{x}}_t,\hat{\bm{y}}_t,\hat{\bm{z}}_t,\hat{\bm{u}}_t,\hat{\bm{e}}_t$. Moreover, we will need the notation
\begin{equation}
\Delta a = a - \tilde{a}
\end{equation}
which means that $\text{var}(a) = \mathbb{E}\left\{\left(\Delta a\right)^2\right\}$.

Based on these notations, we now provide the MO approach for solving \eqref{eq:sc-vector}. The objectives and constraints of MO are listed in Table \ref{tab:mo}, in which we also show the relationship between the original equations in \eqref{eq:sc-vector} and the corresponding equations in MO. The remainder of this section explains the MO method.

\begin{table}[!t]
	\renewcommand{\arraystretch}{1.1}
	%	\linespread{1.5}
	\centering
	\begin{small}
		\caption{Objective and Constraints of MO}
		\vspace{-2mm}
		\begin{tabular}{ccc}
			\hline\hline
			&Original SC \eqref{eq:sc-vector}&MO\\
			\hline
			Decision Variables&$\bm{u}_t^0,\bm{K}$&$\bm{u}_t^0,\bm{K}$  \\
			Objective&\eqref{eq:scp-vector-obj}&\eqref{eq:sc-obj-moment}\\
			Stochastic Resources&\eqref{eq:sc-vector-ito}&\eqref{eq:ito-extend-cov}\\
			EUs&\eqref{eq:eu-vector}&\eqref{eq:eu-1st}\eqref{eq:eu-2nd}\\
			Control Policy&\eqref{eq:sc-control-policy}&\eqref{eq:sc-control-policy-1st}\eqref{eq:sc-control-policy-2nd}\\
			Network Constraints&\eqref{eq:network-vector}\eqref{eq:network-lv-pq-2}&\eqref{eq:network-vector-1st}\eqref{eq:network-lv-pq-1st-approx}\eqref{eq:x-2nd}\\
			Inequality Constraints&\eqref{eq:sc-constraints-x}\eqref{eq:sc-constraints-y}\eqref{eq:sc-constraints-e}\eqref{eq:sc-constraints-p-s}&\eqref{eq:sc-constraints-moment}\\
			\hline \hline
		\end{tabular}\label{tab:mo}
		\vspace{-3mm}
	\end{small}
\end{table}

\vspace{-2mm}
\subsection{Reformulating Objective in \eqref{eq:scp-vector-obj}}
According to the fact that $\bm{R}^U$ is diagonal, we have
\begin{equation}
\mathbb{E}\left\{\bm{u}_t^\top \bm{R}^U\bm{u}_t\right\} = \bm{\tilde{u}}_t^\top \bm{R}^U \bm{\tilde{u}}_t + \bm{\hat{u}}_t^\top \bm{R}^U \bm{\hat{u}}_t
\end{equation}
and the same discussion can be applied to $\mathbb{E}\left\{\bm{e}_T^\top\bm{R}^E\bm{e}_T\right\}$ and $\mathbb{E}\left\{\bm{x}_t^\top\bm{R}^X\bm{x}_t\right\}$. Therefore, we have
\begin{equation} \label{eq:sc-obj-moment}
\begin{split}
J = & \int_{t\in\mathcal{T}}\left(\bm{H}\tilde{\bm{x}}_t+\bm{\tilde{x}}_t^\top \bm{R}^X \bm{\tilde{x}}_t + \bm{\hat{x}}_t^\top \bm{R}^X \bm{\hat{x}}_t\right)dt \\
&+\int_{t\in\mathcal{T}}\left(\bm{\tilde{u}}_t^\top \bm{R}^U \bm{\tilde{u}}_t + \bm{\hat{u}}_t^\top \bm{R}^U \bm{\hat{u}}_t\right)dt \\
&+\bm{\tilde{e}}_T^\top\bm{R}^E\bm{\tilde{e}}_T+\bm{\hat{e}}_T^\top\bm{R}^E\bm{\hat{e}}_T
\end{split}
\end{equation}

\vspace{-2mm}
\subsection{Reformulating Stochastic Resources in \eqref{eq:sc-vector-ito}} \label{section:mo-stochastic}
Here, we must consider the spatial correlation and the temporal correlation, the former of which is described by the covariance matrix, while the latter of which must be considered together with the temporal correlation of EUs in \eqref{eq:eu-vector}. To address the temporal correlation, we define an auxiliary vector as follows:
\begin{equation} \label{eq:eta}
\dot{\bm{\eta}}_t = -\bm{\alpha}\bm{\eta}_t + \bm{\beta}\bm{\xi}_t
\end{equation}
We assume that $\bm{\eta}_0=\bm{0}$. Note that \eqref{eq:eta} is different from \eqref{eq:eu-vector} because it is independent of the decision variable $\bm{u}_t^0$ and $\bm{K}$; however, the next subsection shows that the statistics of $\bm{e}_t$ are determined by the statistics of $\bm{\eta}_t$. 

We can rewrite \eqref{eq:xi-ito-vector} and \eqref{eq:eta} as
\begin{equation} \label{eq:xi-eta-ito}
d
\begin{bmatrix}
\bm{\xi}_t \\
\bm{\eta}_t
\end{bmatrix}
= 
\begin{bmatrix}
\bm{\mu}\left(\bm{\xi}_t\right) \\
-\bm{\alpha}\bm{\eta}_t + \bm{\beta}\bm{\xi}_t
\end{bmatrix}
dt + 
\begin{bmatrix}
\bm{\sigma}\left(\bm{\xi}_t\right) \\
\bm{0}
\end{bmatrix}
d\bm{W}_t
\end{equation}
which is also an It\^{o} process. The statistics needed here include the expectation, defined by $\tilde{\bm{\xi}}_t=\mathbb{E}\bm{\xi}_t$ (note that $\mathbb{E}\bm{\eta}_t=0$), and the covariance matrix, defined by
\begin{equation} \label{eq:ito-extend-cov}
\begin{split}
\mathcal{M}_t &= \mathbb{E}\left\{
\begin{bmatrix}
\bm{\xi}_t - \bm{\tilde{\xi}}_t \\
\bm{\eta}_t - \bm{\tilde{\eta}}_t
\end{bmatrix}
\begin{bmatrix}
\bm{\xi}_t - \bm{\tilde{\xi}}_t \\
\bm{\eta}_t - \bm{\tilde{\eta}}_t
\end{bmatrix}^\top
\right\} 
%&= 
%\begin{bmatrix}
%\mathbb{E}\left(\bm{\xi}_t - \bm{\tilde{\xi}}_t\right)\left(\bm{\xi}_t - \bm{\tilde{\xi}}_t\right)^\top & \mathbb{E}\left(\bm{\xi}_t - \bm{\tilde{\xi}}_t\right)\left(\bm{\eta}_t - \bm{\tilde{\eta}}_t\right)^\top \\
%\mathbb{E}\left(\bm{\eta}_t - \bm{\tilde{\eta}}_t\right)\left(\bm{\xi}_t - \bm{\tilde{\xi}}_t\right)^\top & \mathbb{E}\left(\bm{\eta}_t - \bm{\tilde{\eta}}_t\right)\left(\bm{\eta}_t - \bm{\tilde{\eta}}_t\right)^\top 
%\end{bmatrix} \\
=
\begin{bmatrix}
\mathcal{M}_t^{\xi\xi} &\mathcal{M}_t^{\xi\eta} \\
\mathcal{M}_t^{\eta\xi} & \mathcal{M}_t^{\eta\eta}
\end{bmatrix}
\end{split}
\end{equation}
where $\mathcal{M}_t^{\eta\xi} = \left(\mathcal{M}_t^{\xi\eta}\right)^\top$. In \eqref{eq:ito-extend-cov}, $\mathcal{M}_t^{\xi\xi}$ describes the spatial correlation of $\bm{\xi}_t$, while the other parts describe the temporal correlation that are necessary in MO.

%\subsubsection{Calculating Statistics}
A simple approach to obtaining $\bm{\tilde{\xi}}_t$ and $\mathcal{M}_t$ is the simulation approach, which is time-consuming. However, it is shown in \cite{chen2018stochastic} that the statistics of It\^{o} processes can be efficiently computed by series expansion. Note that $\bm{\tilde{\xi}}_t$ and $\mathcal{M}_t$ are independent of decision variables $\bm{u}_t^0$ and $\bm{K}$; hence, we assume they are given in the following subsections.

\vspace{-2mm}
\subsection{Reformulating Control Policy \eqref{eq:sc-control-policy} and EUs \eqref{eq:eu-vector}}
By taking the first-order moment in \eqref{eq:sc-control-policy}, we have
\begin{equation} \label{eq:sc-control-policy-1st}
\bm{\tilde{u}}_t = \bm{u}_t^0 + \bm{K}\bm{\tilde{\xi}}_t
\end{equation}
\begin{equation}
\Delta \bm{u}_t = \bm{K} \Delta\bm{\xi}_t
\end{equation}
Therefore, we have
\begin{equation}
\mathbb{E}\Delta \bm{u}_t\Delta \bm{u}_t^\top = \bm{K} \mathcal{M}_t^{\xi\xi}\bm{K}^\top
\end{equation}
Moreover, $\hat{\bm{u}}_t$ is the square root of the diagonal of $\mathbb{E}\Delta \bm{u}_t\Delta \bm{u}_t^\top$:
\begin{equation} \label{eq:sc-control-policy-2nd}
\hat{\bm{u}}_t = \sqrt{diag\left\{\bm{K} \mathcal{M}_t^{\xi\xi}\bm{K}^\top\right\}}
\end{equation}

Now, we discuss the moment-form of $\bm{e}_t$. It is easy to show that $\tilde{\bm{e}}_t$ and $\Delta \bm{e}_t$ satisfy
\begin{align}
\dot{\bm{\tilde{e}}} _t &= -\bm{\alpha \tilde{e}}_t + \bm{\beta \tilde{u}}_t \label{eq:eu-1st}\\
\Delta \dot{\bm{e}}_t &= - \bm{\alpha}\Delta\bm{e}_t + \bm{\beta}\bm{K}\Delta\bm{\xi}_t
\end{align}
Then, it is clear that $\Delta \bm{e}_t = \bm{K} \bm{\eta}_t$,\footnote{Such statement needs the assumption that $\bm{\alpha K}=\bm{K\alpha}$ and $\bm{\beta K}=\bm{K\beta}$. However, this requirement is easy to meet if we convert $\bm{\alpha},\bm{\beta},\bm{K}$ to (larger-order) block-diagonal matrices.}; therefore, we have
\begin{equation} \label{eq:eu-2nd}
\hat{\bm{e}}_t= \sqrt{diag\left\{\bm{K} \mathcal{M}_t^{\eta\eta}\bm{K}^\top\right\}}
\end{equation}
%where the ``$\geq$'' is via the same skill as \eqref{eq:sc-control-policy-2nd-soc}.

\vspace{-2mm}
\subsection{Reformulating Network Constraints}
The network constraints include the linear constraints \eqref{eq:network-vector} and the quadratic constraints \eqref{eq:network-lv-pq-2}, the first-order moments of which are
\begin{equation} \label{eq:network-vector-1st}
\bm{\tilde{x}}_t = \bm{A}_y\bm{\tilde{y}}_t +\bm{A}_z\bm{\tilde{z}}_t+\bm{A}_d\bm{d}_t + \bm{A}_e\bm{\tilde{e}}_t + \bm{A}_u\bm{\tilde{u}}_t
\end{equation}
\begin{equation} \label{eq:network-lv-pq-1st}
\tilde{l}_{ij,t}\tilde{v}_{i,t} + \text{cov}(l_{ij,t},v_{i,t})=\tilde{P}_{ij,t}^2+\tilde{Q}_{ij,t}^2+\hat{P}_{ij,t}^2+\hat{Q}_{ij,t}^2,\forall (i,j)\in \mathcal{E}
\end{equation}

To obtain a convex version of \eqref{eq:network-lv-pq-1st}, we claim that 
$\text{cov}(l_{ij,t},v_{i,t})$ can be ignored. Actually, we have
\begin{equation} \label{eq:cov-lv-ignore}
\text{cov}(l_{ij,t},v_{i,t}) \ll \hat{P}_{ij,t}^2+\hat{Q}_{ij,t}^2
\end{equation}
of which the explanation is provided in 
\begin{reduced}
the extended version \cite[Appendix A]{Chen2019Optimal}.
\end{reduced}
\begin{extended}
Appendix \ref{appendix:proof-cov-lv-bound}.
\end{extended}
Thus, \eqref{eq:network-lv-pq-1st} can be replaced by
\begin{equation} \label{eq:network-lv-pq-1st-approx}
\tilde{l}_{ij,t}\tilde{v}_{i,t} =\tilde{P}_{ij,t}^2+\tilde{Q}_{ij,t}^2+\hat{P}_{ij,t}^2+\hat{Q}_{ij,t}^2,\forall (i,j)\in \mathcal{E}
\end{equation}

The second-order moment of $\bm{x}_t$ satisfies
\begin{equation} \label{eq:x-2nd}
\hat{\bm{x}}_t = \sqrt{diag\left\{
	 \bm{L}\mathcal{M}_t\bm{L}^\top\right\}}
\end{equation}
where $\bm{L} = \begin{bmatrix}
\bm{A}_\xi+\bm{A}_u\bm{K} & \bm{A}_e\bm{K}
\end{bmatrix}$. The derivation of \eqref{eq:x-2nd} is provided in 
\begin{reduced}
	the extended version \cite[Appendix B]{Chen2019Optimal}.
\end{reduced}
\begin{extended}
	Appendix \ref{appendix:x-2nd}.
\end{extended}

\vspace{-2mm}
\subsection{Reformulating Inequality Constraints}
It is shown in \cite{calafiore2006distributionally} that a chance constraint can be approximately described by a second-order cone constraint. In this approach, the constraints in \eqref{eq:sc-constraints-x}$\sim$\eqref{eq:sc-constraints-p-s} need to be handled row by row. For simplicity, we assume that $a$ is a variable and take $a \leq \bar{a}$ as an example.

According to \cite{calafiore2006distributionally}, the second-order-cone formulation of the constraint is
\begin{equation} \label{eq:sc-constraint-soc}
\tilde{a} + \kappa_\gamma \hat{a} \leq \bar{a}
\end{equation}

Therefore, the inequality constraints can be transformed into
\begin{equation} \label{eq:sc-constraints-moment}
\begin{array}{c}
\bm{\underline{x}} + \kappa_\gamma\bm{\hat{x}}_t \leq \bm{\tilde{x}}_t \leq \bm{\bar{x}} - \kappa_\gamma\bm{\hat{x}}_t \\
\bm{\underline{y}} + \kappa_\gamma\bm{\hat{y}}_t \leq \bm{\tilde{y}}_t \leq \bm{\bar{y}} - \kappa_\gamma\bm{\hat{y}}_t \\
\bm{\underline{e}} + \kappa_\gamma\bm{\hat{e}}_t \leq \bm{\tilde{e}}_t \leq \bm{\bar{e}} - \kappa_\gamma\bm{\hat{e}}_t \\
\bm{C}_u\tilde{\bm{u}}_t + \bm{C}_\xi\tilde{\bm{\xi}}_t + \kappa_\gamma\left(\left|\bm{C}_u\right|\hat{\bm{u}}_t+\left|\bm{C}_\xi\right|\hat{\bm{\xi}}_t\right) \leq \bm{D}
\end{array}
\end{equation}

\vspace{-2mm}
\subsection{Summary}
The MO approach can be summarized as
\begin{equation} \label{eq:mo}
\begin{split}
\text{Objective:~}&\eqref{eq:sc-obj-moment}\\
\text{Constraints:~}&\eqref{eq:ito-extend-cov}\eqref{eq:sc-control-policy-1st}\eqref{eq:sc-control-policy-2nd}\eqref{eq:eu-1st}\eqref{eq:eu-2nd}\eqref{eq:network-vector-1st}\eqref{eq:network-lv-pq-1st-approx}\eqref{eq:x-2nd}\eqref{eq:sc-constraints-moment}
\end{split}
\end{equation}
Various nonconvex constraints exist, i.e., \eqref{eq:sc-control-policy-2nd}\eqref{eq:eu-2nd}\eqref{eq:network-lv-pq-1st-approx}\eqref{eq:x-2nd}. However, it is easy to obtain their exact convex relaxations 
\begin{reduced}
	(see Appendix C in the extended version \cite{Chen2019Optimal}).
\end{reduced}
\begin{extended}
(see Appendix \ref{appendix:convexity}).
\end{extended}
Therefore, \eqref{eq:mo} is a convex optimization problem.

Now, we discuss the computational burden of MO. It is clear that for a certain variable, say $a$, in \eqref{eq:sc-vector}, there are two corresponding variables in MO, i.e., $\tilde{a}$ and $\hat{a}$ respectively. Therefore, the number of variables of MO is approximately twice that of the original SC problem. Moreover, each constraint in \eqref{eq:sc-vector} corresponds with one or two constraints in MO, regarding the first-order and second-order moments. Therefore, the number of constraints of MO is less than twice that of the original SC problem. In contrast, the scale of traditional SPBC algorithms is proportional to the number of scenarios, which is usually large for an accurate estimation of the SC problem with spatial and temporal correlation. In summary, the MO approach reduces the computational burden of SC to be comparable with the corresponding deterministic control problem.

\section{Case Study} \label{section:case}
This section provides a test case in an IEEE 123-bus distribution network \cite{ieee123}. We evaluate the optimal control scheme provided by the MO approach and then discuss the impacts of spatial and temporal correlation. Moreover, the comparison between the proposed MO approach and several existing approaches shows the effectiveness and efficiency of the MO approach.

\subsection{Case Settings}
\begin{figure}
	\centering
	\includegraphics[width=\columnwidth]{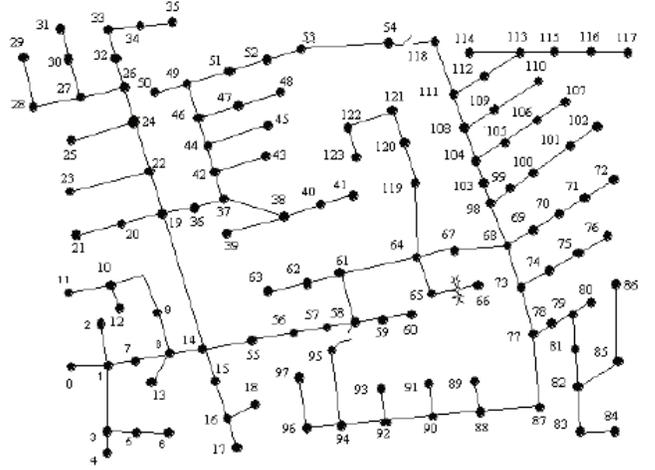}
	\caption{IEEE 123-bus system.}
	\label{fig:ieee-123}
\end{figure}

We consider the IEEE 123-bus system, as shown in Fig. \ref{fig:ieee-123}. The parameters of the IEEE 123-bus system can be found in \cite{ieee123}. The nominal capacity of the system is 10~MVA, and the nominal voltage is 10~kV. We assume that the voltage limit of each bus is $10\pm 0.5$~kV. The stochastic resources are 3 wind generators on Buses 11, 62, and 66, the capacity of each of which is 20~MVA, and 3 PV generators on Buses 72, 75, and 114, the capacity of each of which is 10~MVA. We use \eqref{eq:ito-example-3} to describe the stochastic resources, where $\bm{\mu}$ is an affine function, and $\bm{\sigma}$ is constant.  We assume that the predicted values are obtained by persistent prediction in 1 h \cite{bludszuweit2008statistical}, and the parameters are obtained via the parameter estimation method provided in \cite{du2017parameter}. We assume that there exists a correlation between Buses 62 and 66 and a correlation between Buses 72 and 75, of which the correlation coefficients are both 0.5. The parameters are provided in
\begin{reduced}
	the extended version \cite[Appendix D]{Chen2019Optimal}.
\end{reduced}
\begin{extended}
	Appendix \ref{appendix:parameter}.
\end{extended}
Moreover, there is an EU at Bus 62 (5 MW $\times$ 4 h).

The objective of MO is as shown in \eqref{eq:scp-obj}, where the price is \$1/kWh from 08:00$\sim$20:00, and \$0.5/kWh during the rest of the day. $\bm{R}$ is a diagonal matrix whose diagonal elements are all 1, and $R^V=1, R^E=0.1$.
%In order to ensure feasibility, we use soft-constraint technique to deal with the inequality constraints in \eqref{eq:sc-constraints-moment} \cite{valverde2013model}, in which the penalty factor is 100.
The time step is 15 minutes; the control horizon is 1 day; and the $\bm{u}_t^0$ and $\bm{K}$ are updated every 4 hours in order to maintain a good performance. 

\subsection{Simulation Results}
To evaluate the control scheme obtained by the MO approach, we use a Monte Carlo simulation with 1000 scenarios to calculate the objective function under the control scheme. The objective value under the optimal control scheme is $J=249$ k\$. In contrast, if we let $\bm{K}=\bm{0}$ and only consider $\bm{u}_t^0$, the result is $J'=$ 272 k\$. In fact, $J'$ is the objective under deterministic control schemes, which means that the controlled units do not respond to any disturbances of stochastic resources. The results show that stochastic control scheme performs better than deterministic control schemes.

The value of $\bm{K}$ shows the relationship between $\bm{u}_t$ and $\bm{\xi}_t$. In this case, we have

\begin{equation*}
\begin{small}
\bm{K} = 
\begin{bmatrix}
%0 & 0 & 0 & 0 & 0 & 0 \\
-0.068 & -0.087 & -0.096 & -0.001 & -0.002 & -0.001 \\
-0.087 & -0.184 & -0.195 & -0.003 & -0.002 & -0.003 \\
-0.096 & -0.916 & -0.421 & -0.003 & -0.004 & -0.002 \\
-0.001 & -0.003 & -0.002 & -0.085 & -0.069 & -0.024 \\
-0.002 & -0.002 & -0.004 & -0.069 & -0.096 & -0.025 \\
-0.001 & -0.003 & -0.002 & -0.024 & -0.025 & -0.043 \\
-0.102 & -0.203 & -0.184 & -0.014 & -0.017 & -0.006
\end{bmatrix}
\end{small}
\end{equation*}
where the order of the columns is the output of stochastic resources at Buses 11, 62, 66, 72, 75, and 114, and the order of the rows is the reactive power of stochastic resources and the output of the EU. It is clear that the control scheme is a negative feedback control scheme. Moreover, the values of $\bm{K}$ show the correlation between these variables, and units at closer buses share larger coefficients. For example, the Row 7, Column 2 of $\bm{K}$ describing the sensitivity of the EU output with respect to the DG at Bus 62, is relatively larger.

\begin{figure}
	\centering
	\includegraphics[width=\columnwidth]{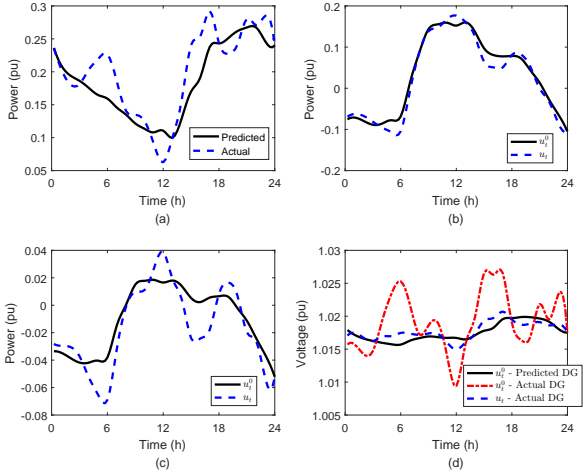}
	\caption{Simulation results at Bus 62. (a) power of wind generations; (b) reactive power of wind generations; (c) power of EU; (d) voltage of Bus 62.}
	\label{fig:results}
\end{figure}

Fig. \ref{fig:results} shows the curves at Bus 62 in a certain scenario. The negative feedback control scheme is shown by Fig. \ref{fig:results}(a)(b)(c), where lower wind power leads to larger control output. Fig.~\ref{fig:results}(d) shows the effect of the feedback control, where the black curve is the voltage profile under perfect prediction; the red curve is the voltage profile under uncertainty, but with feedback coefficient $\bm{K}=\bm{0}$; and the blue curve is the voltage profile under the optimal feedback control. It is clear that feedback control improves the voltage profile at Bus 62.

\subsection{Impacts of Correlation} \label{section:case-correlation}
\begin{figure}
	\centering
	\includegraphics[width=\columnwidth]{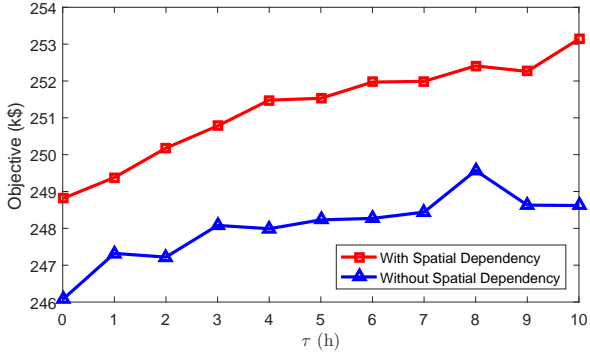}
	\caption{Impacts of correlation.}
	\label{fig:correlation}
\end{figure}

Here, we discuss the relationship between the correlation and the control performance. For the spatial correlation, we consider $\bm{\sigma}'=diag\left\{2.98, 7.52, 4.51, 1.42, 3.75, 2.76\right\}$. It is clear that the diffusion coefficients $\bm{\sigma}$ and $\bm{\sigma}'$ result in the same variance of $\bm{\xi}_t$, but the stochastic resources under $\bm{\sigma}'$ are spatially independent. For temporal correlation, Example \ref{eg:ito-1} shows that it will not change anything except the temporal correlation if we multiply $\bm{\mu}$ by $1/\tau$ and $\bm{\sigma}$ by $1/\sqrt{\tau}$ simultaneously; therefore, $\tau$ can be regarded as a measure of temporal correlation ($\tau=1$ for default situation). Fig. \ref{fig:correlation} shows the objective under different spatial and temporal correlation. It is shown that larger (spatial or temporal) correlation leads to worse control performance, and we now explain this result. Some of the variables, such as $\bm{e}_t$, are related to the integration or sum of elements in $\bm{\xi}_t$. However, the uncertainty of the sum of stochastic variables is influenced not only by the variance of each variable but also by the correlation. Moreover, the correlation in this case is positive and hence will lead to larger uncertainty and a larger objective value.

\begin{figure}
	\centering
	\includegraphics[width=\columnwidth]{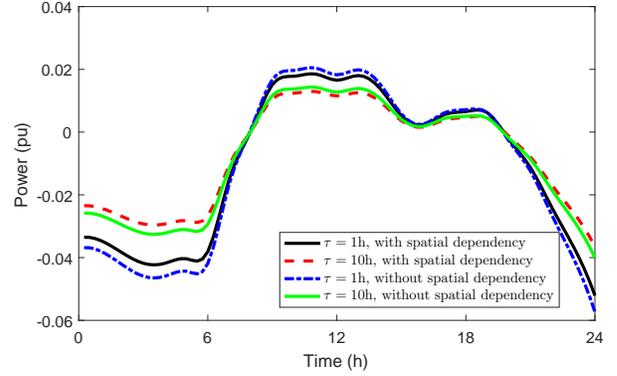}
	\caption{Control schemes under different correlation.}
	\label{fig:correlation-example}
\end{figure}

Fig. \ref{fig:correlation-example} shows the $\bm{u}_t^0$ (taking the EU as an example) under different spatial and temporal correlations. It is also clearly shown that larger spatial and temporal correlation results in smaller controller output, since the controller must reserve more capacity for uncertainty. The impacts of spatial and temporal correlation can be considered in the proposed approach.

\subsection{Comparisons with Existing Approaches}
This section compares MO with other existing approaches. Here, we use the deterministic control (DC), MPC, and SPBC as benchmarks. In the DC approach, the control scheme is obtained based on the predicted values of stochastic resources, and the uncertainties are ignored. MPC performs DC in a receding-horizon manner, and update the prediction value at each time step. The prediction horizon of MPC is 4 hours. We perform the SPBC approach with 100 scenarios and 1000 scenarios, denoted by SPBC(20) and SPBC(100). These scenarios are obtained by the scenario-reduction method provided in \cite{li2016scenario}. We compare the performance and the computational burden of these control methods.

\subsubsection{Performance of Different Methods}
Table \ref{tab:comparison} shows the performance of different methods, demonstrating that MO and SPBC(100) perform best and that DC performs the worst. DC does not consider the uncertainty of the stochastic resources and hence achieves the worst performance, which is also supported by Fig. \ref{fig:results}. MO and SPBC both consider the uncertainty explicitly and perform well. Moreover, SPBC(20) does not perform as well as MO because 20 scenarios are too few to describe the correlation of the stochastic resources. Although MPC does not explicitly consider the uncertainty, the receding-horizon manner improves its performance. Nevertheless, MPC does not perform as well as MO.

\subsubsection{Computational Burden}
Table \ref{tab:comparison} shows the computational burden of these methods. Since MPC is a receding-horizon control method, while the others are not, we use the per-step computational time for a fair comparison. It is clear that DC is the fastest because it does not consider the uncertainty. The computational time of MO is about twice that of DC, significantly smaller than MPC, SPBC(20) and SPBC(100), which shows the advantage of the proposed method over existing methods. Specifically, since MO and SPBC(100) achieve similar performance, it can be concluded that MO reduces the computational time by 99.5\% without sacrificing performance.

In summary, the proposed MO approach achieves a good trade-off between the control performance and the computational burden. In contrast, DC is computationally efficient but performs worse, while MPC/SPBC performs well but incurs an extremely large computational burden. Therefore, the proposed MO significantly outperforms the existing methods and has attractive potential in the control of DERs under uncertainty. 

\begin{table}[!t] 
	%	\vspace{-5mm}
	\renewcommand{\arraystretch}{1.1}
	%	\linespread{1.5}
	\centering
	\begin{small}
		\caption{Benchmarks}
		\begin{tabular}{c|cc}
			\hline\hline
			\multirow{2}{*}{Approach}& Computation&Objective (k\$)\\
			&Time (min per step)&Value\\
			\hline
			MO&0.5&249.5\\
			\hline
			DC&0.2&276.3\\
			\hline
			MPC&12&257.6\\
			\hline
			SPBC(20)&25&262.4\\
			\hline
			SPBC(100)&107&250.8\\
			\hline \hline
		\end{tabular}
		\label{tab:comparison}
	\end{small}
	%	\vspace{-4mm}
\end{table}

\section{Conclusion} \label{section:conclusion}
This paper presents an MO approach for the efficient control of DERs in distribution networks. We first model the stochastic resources by It\^{o} processes, which describe the spatial and temporal correlation of the stochastic resources. The It\^{o} process model is also in the same form as the characteristics of EUs; hence, the temporal correlation of the stochastic resources and EUs can be considered in a unified way. Based on the covariance matrix obtained by the It\^{o} process model, we transform the SC problem into a deterministic optimization problem with respect to the first-order and second-order moments of the system variables, whose scale is approximately twice that of the corresponding deterministic control problem. The proposed MO approach solves the SC problem in a computationally efficient way and outperforms existing approaches such as DC, MPC, and SPBC.

\begin{extended}

\appendices
\section{Explanation of Inequality \eqref{eq:cov-lv-ignore}} \label{appendix:proof-cov-lv-bound}
	Consider the first-order deviation of \eqref{eq:network-lv-pq-2}, and ignore higher-order deviations:
	\begin{equation}
	\tilde{l}_{ij,t}\Delta v_{i,t} + \tilde{v}_{i,t}\Delta l_{ij,t} = 2\tilde{P}_{ij,t}\Delta P_{ij,t} + 2\tilde{Q}_{ij,t}\Delta Q_{ij,t}
	\end{equation}
	Therefore,
	\begin{equation}
	\begin{split}
	\tilde{v}_{i,t}&\Delta l_{ij,t}\Delta v_{i,t} = 2\tilde{P}_{ij,t}\Delta P_{ij,t}\Delta v_{i,t} \\
	&+ 2\tilde{Q}_{ij,t}\Delta Q_{ij,t}\Delta v_{i,t}- \tilde{l}_{ij,t}\left(\Delta v_{i,t}\right)^2
	\end{split}
	\end{equation}
	By taking expectations in both sides, we have
	\begin{equation}
	\begin{split}
	\text{cov}&\left(l_{ij,t},v_{i,t}\right) = 2\frac{\tilde{P}_{ij,t}}{\tilde{v}_{i,t}}\text{cov}\left( P_{ij,t}, v_{i,t}\right) \\
	&+ 2\frac{\tilde{Q}_{ij,t}}{\tilde{v}_{i,t}}\text{cov}\left( Q_{ij,t},v_{i,t}\right) - \frac{\tilde{l}_{ij,t}}{\tilde{v}_{i,t}}\hat{v}_{i,t}^2
	\end{split}
	\end{equation}
	By applying Cauchy inequality and using the fact that $\tilde{l}_{ij,t}/\tilde{v}_{i,t} > 0$, we have
\begin{equation} \label{eq:cov-lv-bound}
\text{cov}(l_{ij,t},v_{i,t}) \leq 2 \frac{\hat{v}_{i,t}/\tilde{v}_{i,t}}{\hat{P}_{ij,t}/\left|\tilde{P}_{ij,t}\right|}\hat{P}_{ij,t}^2 + 2 \frac{\hat{v}_{i,t}/\tilde{v}_{i,t}}{\hat{Q}_{ij,t}/\left|\tilde{Q}_{ij,t}\right|}\hat{Q}_{ij,t}^2
\end{equation}

The term $\frac{\hat{v}_{i,t}/\tilde{v}_{i,t}}{\hat{P}_{ij,t}/\left|\tilde{P}_{ij,t}\right|}$ can be interpreted as the sensitivity of the voltage deviation $v_{i,t}$ under the power deviation $P_{ij,t}$, which, in practice, is very small since the relative deviation of the bus voltage is far less than that of the branch power. Therefore, we have $\text{cov}(l_{ij,t},v_{i,t}) \ll \hat{P}_{ij,t}^2+\hat{Q}_{ij,t}^2$.

\section{Second-Order Moment of Network Constraints} \label{appendix:x-2nd}
A major challenge to obtain the second-order-moment form of network constraints is how to avoid higher-order moments. To achieve this, we consider a typical approximation form of the distFlow model \cite{liu2017decentralized}:
\begin{equation} \label{eq:distflow-approx}
\begin{array}{c}
p_{j,t}=P_{jk,t} +\sum_{i:i\to j} \left(P_{ij,t}\right) + g_jv_{j,t},\forall j \\
q_{j,t}=Q_{jk,t} -\sum_{i:i\to j} \left(Q_{ij,t}\right) + b_jv_{j,t},\forall j \\
v_{j,t}= v_{j,t} - 2\left(r_{ij}P_{ij,t}+x_{ij}Q_{ij,t}\right),\forall (i,j)\in \mathcal{E}
\end{array}
\end{equation}
This model assumes negligible line losses and almost flat voltage, and its accuracy has been verified by several recent work \cite{li2016convex,sulc2014optimal}. 

%\begin{remark}
	It must be emphasized that the approximated model \eqref{eq:distflow-approx} is only used to estimate the second-order moments, while the first-order moments in \eqref{eq:network-vector-1st} and \eqref{eq:network-lv-pq-1st-approx} are computed by the exact model shown in \eqref{eq:network-vector} and \eqref{eq:network-lv-pq-2}. In other words, we use an accurate model to estimate the expectations and an approximate model to estimate the errors. Since expectations are usually more important in the objective of stochastic optimization problems, the use of the approximation will have acceptable impacts on the accuracy of the model.
%\end{remark}

It is clear that the vector form of \eqref{eq:distflow-approx} is
\begin{equation}
\bm{x}_t = \bm{A}_\xi\bm{\xi}_t+\bm{A}_d\bm{d}_t + \bm{A}_e\bm{e}_t + \bm{A}_u\bm{u}_t
\end{equation}
By replacing $\bm{e}_t$ and $\bm{u}_t$ by $\bm{\xi}_t$ and $\bm{\eta}_t$, we have
\begin{equation}
\Delta \bm{x}_t = \left(\bm{A}_\xi+\bm{A}_u\bm{K}\right)\Delta\bm{\xi}_t+ \bm{A}_e\bm{K}\bm{\eta}_t
\end{equation}
then \eqref{eq:x-2nd} clearly follows.
\section{Exact Convex Relaxation of MO} \label{appendix:convexity}
The exact convex relaxations of \eqref{eq:sc-control-policy-2nd}\eqref{eq:eu-2nd}\eqref{eq:x-2nd} are
\begin{equation}  \label{eq:conic-relaxation-diag}
\begin{split}
\hat{\bm{u}}_t & \geq \sqrt{diag\left\{\bm{K} \mathcal{M}_t^{\xi\xi}\bm{K}^\top\right\}} \\
\hat{\bm{e}}_t &\geq \sqrt{diag\left\{\bm{K} \mathcal{M}_t^{\eta\eta}\bm{K}^\top\right\}} \\
\hat{\bm{x}}_t &\geq \sqrt{diag\left\{
	\bm{L}\mathcal{M}_t\bm{L}^\top\right\}}
\end{split}
\end{equation}
Here we only explain the first equation as an example. Since $\mathcal{M}_t^{\xi\xi}$ is symmetric, we assume $\mathcal{M}_t^{\xi\xi}=\mathcal{N}_t\mathcal{N}_t^\top$. Considering the $i$-th coordinate of $\bm{u}_t$, denoted by $\bm{u}_{t,i}$, we have 
\begin{equation}
\bm{u}_{t,i}=\sqrt{\left\{\bm{K}_i \mathcal{M}_t^{\xi\xi}\bm{K}_i^\top\right\}}=\left\|\bm{K}_i\mathcal{N}_t\right\|_2
\end{equation}
where $\bm{K}_i$ is the $i$-th row of $\bm{K}$. Clearly, this equation can be exactly relaxed as
\begin{equation}
\bm{u}_{t,i}\geq\sqrt{\left\{\bm{K}_i \mathcal{M}_t^{\xi\xi}\bm{K}_i^\top\right\}}=\left\|\bm{K}_i\mathcal{N}_t\right\|_2
\end{equation}
which is a second-order-cone constraint. And \eqref{eq:conic-relaxation-diag} can be obtained in a similar way.

Moreover, the exact convex relaxation of \eqref{eq:network-lv-pq-1st-approx} is
\begin{equation} \label{eq:conic-relaxation-1st}
\tilde{l}_{ij,t}+\tilde{v}_{i,t} \geq\left\|
\begin{matrix}
2\tilde{P}_{ij,t} \\
2\tilde{Q}_{ij,t} \\
2\hat{P}_{ij,t} \\
2\hat{Q}_{ij,t} \\
\tilde{l}_{ij,t}-\tilde{v}_{i,t}
\end{matrix}
\right\|_2,\forall (i,j)\in \mathcal{E}
\end{equation}
where $\left\|\cdot\right\|_2$ is the operator of 2-norm. This relaxation technique is widely used in the convex relaxation of distFlow \cite{farivar2013branch}.

\section{Parameters of Stochastic Resources} \label{appendix:parameter}
Here we provide the $\bm{\sigma}$, which will influence the correlation of $\bm{\xi}_t$, and will be discussed in Section \ref{section:case-correlation}. Other parameters are omitted due to space constraints.
\begin{equation}
\bm{\sigma} = 
\begin{bmatrix}
2.98 & 0 & 0 & 0 & 0 & 0 \\
0 & 7.52 & 0 & 0 & 0 & 0 \\
0 & 2.25 & 3.91 & 0 & 0 & 0 \\
0 & 0 & 0 & 1.42 & 0 & 0 \\
0 & 0 & 0 & 0 & 3.75 & 0 \\
0 & 0 & 0 & 0 & 1.46 & 2.35
\end{bmatrix}
\end{equation}

\end{extended}

\bibliographystyle{ieeetr}
\bibliography{scp-adn}

\begin{thebibliography}{10}

\bibitem{maycock2014future}
P.~Maycock, ``{The Future of Energy}.'' \url{http://bnef.com/
  Presentations/download/136}, 2014.
\newblock [Online; accessed 16-Aug-2018].

\bibitem{hidalgo2010review}
R.~Hidalgo, C.~Abbey, and G.~Joós, ``A review of active distribution networks
  enabling technologies,'' in {\em IEEE PES General Meeting}, pp.~1--9, July
  2010.

\bibitem{olivier2016active}
F.~Olivier, P.~Aristidou, D.~Ernst, and T.~V. Cutsem, ``Active management of
  low-voltage networks for mitigating overvoltages due to photovoltaic units,''
  {\em IEEE Trans. Smart Grid}, vol.~7, pp.~926--936, March 2016.

\bibitem{li2016convex}
Q.~Li, R.~Ayyanar, and V.~Vittal, ``Convex optimization for des planning and
  operation in radial distribution systems with high penetration of
  photovoltaic resources,'' {\em IEEE Trans. Sustain. Energy}, vol.~7,
  pp.~985--995, July 2016.

\bibitem{sugihara2013economic}
H.~Sugihara, K.~Yokoyama, O.~Saeki, K.~Tsuji, and T.~Funaki, ``Economic and
  efficient voltage management using customer-owned energy storage systems in a
  distribution network with high penetration of photovoltaic systems,'' {\em
  IEEE Trans. Power Syst.}, vol.~28, pp.~102--111, Feb 2013.

\bibitem{bludszuweit2008statistical}
H.~Bludszuweit, J.~A. Dominguez-Navarro, and A.~Llombart, ``Statistical
  analysis of wind power forecast error,'' {\em IEEE Trans. Power Syst.},
  vol.~23, pp.~983--991, Aug 2008.

\bibitem{li2015data}
P.~Li, R.~Dargaville, F.~Liu, J.~Xia, and Y.~Song, ``Data-based statistical
  property analyzing and storage sizing for hybrid renewable energy systems,''
  {\em IEEE Trans. Ind. Electron.}, vol.~62, pp.~6996--7008, Nov 2015.

\bibitem{teng2013optimal}
J.~Teng, S.~Luan, D.~Lee, and Y.~Huang, ``Optimal charging/discharging
  scheduling of battery storage systems for distribution systems interconnected
  with sizeable pv generation systems,'' {\em IEEE Trans. Power Syst.},
  vol.~28, pp.~1425--1433, May 2013.

\bibitem{atwa2010optimal}
Y.~M. Atwa, E.~F. El-Saadany, M.~M.~A. Salama, and R.~Seethapathy, ``Optimal
  renewable resources mix for distribution system energy loss minimization,''
  {\em IEEE Trans. Power Syst.}, vol.~25, pp.~360--370, Feb 2010.

\bibitem{tabone2015modeling}
M.~D. Tabone and D.~S. Callaway, ``Modeling variability and uncertainty of
  photovoltaic generation: A hidden state spatial statistical approach,'' {\em
  IEEE Trans. Power Syst.}, vol.~30, pp.~2965--2973, Nov 2015.

\bibitem{lorenz2009irradiance}
E.~Lorenz, J.~Hurka, D.~Heinemann, and H.~G. Beyer, ``Irradiance forecasting
  for the power prediction of grid-connected photovoltaic systems,'' {\em IEEE
  J. Sel. Topics Appl. Earth Observ.}, vol.~2, pp.~2--10, March 2009.

\bibitem{chen2018unified}
X.~Chen, J.~Lin, C.~Wan, Y.~Song, and H.~Luo, ``A unified frequency-domain
  model for automatic generation control assessment under wind power
  uncertainty,'' {\em IEEE Trans. Smart Grid}, Early Access.

\bibitem{anese2015optimal}
E.~Dall’Anese, S.~V. Dhople, B.~B. Johnson, and G.~B. Giannakis, ``Optimal
  dispatch of residential photovoltaic inverters under forecasting
  uncertainties,'' {\em IEEE Journal of Photovoltaics}, vol.~5, pp.~350--359,
  Jan 2015.

\bibitem{li2018coordinated}
J.~Li, Z.~Xu, J.~Zhao, and C.~Wan, ``A coordinated dispatch model for
  distribution network considering pv ramp,'' {\em IEEE Trans. Power Syst.},
  vol.~33, pp.~1107--1109, Jan 2018.

\bibitem{valverde2013model}
G.~Valverde and T.~V. Cutsem, ``Model predictive control of voltages in active
  distribution networks,'' {\em IEEE Trans. Smart Grid}, vol.~4,
  pp.~2152--2161, Dec 2013.

\bibitem{xing2017model}
X.~Xing, J.~Lin, C.~Wan, and Y.~Song, ``Model predictive control of lpc-looped
  active distribution network with high penetration of distributed
  generation,'' {\em IEEE Trans. Sustain. Energy}, vol.~8, pp.~1051--1063, July
  2017.

\bibitem{meng2015cooperation}
K.~Meng, Z.~Y. Dong, Z.~Xu, and S.~R. Weller, ``Cooperation-driven distributed
  model predictive control for energy storage systems,'' {\em IEEE Trans. Smart
  Grid}, vol.~6, pp.~2583--2585, Nov 2015.

\bibitem{agalgaonkar2015stochastic}
Y.~P. Agalgaonkar, B.~C. Pal, and R.~A. Jabr, ``Stochastic distribution system
  operation considering voltage regulation risks in the presence of pv
  generation,'' {\em IEEE Trans. Sustain. Energy}, vol.~6, pp.~1315--1324, Oct
  2015.

\bibitem{mesbah2016stochastic}
A.~Mesbah, ``{Stochastic Model Predictive Control : An Overview and
  Perspectives for Future Research},'' {\em IEEE Commun. Mag.}, no.~December,
  2016.

\bibitem{fu2016multiobjective}
Y.~Fu, M.~Liu, and L.~Li, ``Multiobjective stochastic economic dispatch with
  variable wind generation using scenario-based decomposition and asynchronous
  block iteration,'' {\em IEEE Trans. Sustain Energy}, vol.~7, pp.~139--149,
  Jan 2016.

\bibitem{zhu2014decomposition}
D.~Zhu and G.~Hug-Glanzmann, ``Decomposition methods for stochastic optimal
  coordination of energy storage and generation,'' in {\em 2014 IEEE PES
  General Meeting | Conference Exposition}, pp.~1--5, July 2014.

\bibitem{lamberton2011introduction}
D.~Lamberton and B.~Lapeyre, {\em Introduction to stochastic calculus applied
  to finance}.
\newblock Chapman and Hall/CRC, 2011.

\bibitem{chen2018stochastic}
X.~Chen, J.~Lin, F.~Liu, and Y.~Song, ``Stochastic assessment of agc systems
  under non-gaussian uncertainty,'' {\em IEEE Trans. Power Syst.}, Early
  Access.

\bibitem{du2017parameter}
X.-L. Du, J.-G. Lin, and X.-Q. Zhou, ``Parameter estimation for multivariate
  diffusion processes with the time inhomogeneously positive semidefinite
  diffusion matrix,'' {\em Communications in Statistics-Theory and Methods},
  vol.~46, no.~22, pp.~11010--11025, 2017.

\bibitem{ait2002maximum}
Y.~A{\"\i}t-Sahalia, ``Maximum likelihood estimation of discretely sampled
  diffusions: a closed-form approximation approach,'' {\em Econometrica},
  vol.~70, no.~1, pp.~223--262, 2002.

\bibitem{ait2008closed}
Y.~Ait-Sahalia {\em et~al.}, ``Closed-form likelihood expansions for
  multivariate diffusions,'' {\em The Annals of Statistics}, vol.~36, no.~2,
  pp.~906--937, 2008.

\bibitem{yoshida1992estimation}
N.~Yoshida, ``Estimation for diffusion processes from discrete observation,''
  {\em Journal of Multivariate Analysis}, vol.~41, no.~2, pp.~220--242, 1992.

\bibitem{farivar2013branch}
M.~Farivar and S.~H. Low, ``Branch flow model: Relaxations and
  convexification—part i,'' {\em IEEE Trans. Power Syst.}, vol.~28,
  pp.~2554--2564, Aug 2013.

\bibitem{skaf2010design}
J.~Skaf and S.~P. Boyd, ``{Design of affine controllers via convex
  optimization},'' {\em IEEE Trans. Autom. Control}, vol.~55, no.~11,
  pp.~2476--2487, 2010.

\bibitem{calafiore2006distributionally}
G.~C. Calafiore and L.~El~Ghaoui, ``On distributionally robust
  chance-constrained linear programs,'' {\em Journal of Optimization Theory and
  Applications}, vol.~130, no.~1, pp.~1--22, 2006.

\bibitem{ieee123}
``{IEEE 123-Bus System}.''
  \url{http://sites.ieee.org/pes-testfeeders/resources/}.

\bibitem{li2016scenario}
J.~Li, F.~Lan, and H.~Wei, ``A scenario optimal reduction method for wind power
  time series,'' {\em IEEE Trans. Power Systems}, vol.~31, pp.~1657--1658,
  March 2016.

\bibitem{liu2017decentralized}
H.~J. Liu, W.~Shi, and H.~Zhu, ``Decentralized dynamic optimization for power
  network voltage control,'' {\em IEEE Trans. Signal Inf. Process. Netw.},
  vol.~3, pp.~568--579, Sept 2017.

\bibitem{sulc2014optimal}
P.~Šulc, S.~Backhaus, and M.~Chertkov, ``Optimal distributed control of
  reactive power via the alternating direction method of multipliers,'' {\em
  IEEE Trans. Energy Conver.}, vol.~29, pp.~968--977, Dec 2014.

\end{thebibliography}

\end{document}